\documentclass[english]{cccconf}
\usepackage[comma,numbers,square,sort&compress]{natbib}
\usepackage{bm,algorithmicx,algpseudocode,threeparttable,slashbox,booktabs}
\usepackage{graphicx,times,amsmath}
\usepackage{stmaryrd}
\usepackage{amsfonts}
\usepackage{hyperref}

\usepackage{fancyvrb}
\usepackage{color}

\begin{document}

\title{A matlab toolbox  for continuous state transition algorithm}
\author{Xiaojun Zhou, Chunhua Yang, Weihua Gui}



\affiliation{School of Information Science and Engineering,
        Cental South University, Changsha 410083, P.~R.~China
        \email{michael.x.zhou@csu.edu.cn}}

\maketitle

\begin{abstract}
State transition algorithm (STA) has been emerging as a novel stochastic method for global optimization in recent few years. 
To make better understanding of continuous STA, a matlab toolbox for continuous STA has been developed. 
Firstly, the basic principles of continuous STA are briefly described.
Then, a matlab implementation of the standard continuous STA is explained, with 
several instances given to show how to use to the matlab toolbox to minimize an optimization problem with bound constraints.
In the same while, a link is provided to download the matlab toolbox via available resources. 
\end{abstract}


\footnotetext{This work is supported by National Natural Science
Foundation (NNSF) of China under Grant No. 61503416, 61533020, 61590921.}

\section{Introduction}
State transition algorithm (STA) \cite{Zhou2011initial, Zhou2011new, Zhou2012, Zhou2013, Yang2013, Zhou2015discrete} has been emerging as a novel stochastic method for global optimization in recent few years,
and it has found applications in nonlinear system identification, optimized controller design, water distribution networks, energy conservation optimization, optimized controller design, signal processing, etc \cite{Zhou2014,Zhou2015optimal,Wang2016optimal,Wang2016sta}.
In state transition algorithm, a solution to
an optimization problem is considered as a state, and an update of a solution can be
regarded as a state transition.
Unlike other population-based stochastic optimization techniques, such as
genetic algorithm, particle swarm optimization, differential evolution, etc, the basic state transition algorithm is an individual-based optimization method.
Based on an incumbent best solution, a neighborhood with special property will be formed automatically
when using certain state transformation operator.
A variety of state transformation operators, for example,
rotation, translation, expansion, and axesion in continuous STA, or swap, shift, symmetry and
substitute in discrete STA,
are designed purposely for both global and local search. On the basis of the neighborhood,
then, a sampling technique is used to generate a candidate set, and
the next best solution is updated by using a selection technique based on previous best solution and
the candidate set. This process is repeated until some terminal conditions are satisfied.

In this paper, we focus on the continuous STA for the following global optimization problem
\begin{eqnarray}
\min_{\bm x \in \Omega} f(\bm x)
\end{eqnarray}
where $\bm x \in \mathbb{R}^n$, $\Omega \subseteq \mathbb{R}^n$ is a closed and compact set, which is usually composed of lower and upper bounds of $\bm x$.

To make better understanding of continuous STA, a matlab toolbox for continuous STA for global optimization problem with 
bound constraints has been developed.
In the rest of this paper, we will describe the continuous STA and its corresponding matlab toolbox in detail. 

\section{The basic principles of continuous STA}
By referring to state space representation, on the basis of current state $\bm x_k$, the unified form of generation of a new state $\bm x_{k+1}$ in state transition algorithm can be described as follows:
\begin{equation}
\left \{ \begin{array}{ll}
\bm x_{k+1}= A_{k} \bm x_{k} + B_{k} \bm u_{k}\\
y_{k+1}= f(\bm x_{k+1})
\end{array} \right.,
\end{equation}
where $\bm x_{k} = [x_1, x_2, \cdots, x_n]^T$ stands for a state, corresponding to a solution of an optimization problem; $\bm u_{k}$ is a function of $\bm x_{k}$ and historical states;
$A_{k}$ and
$B_{k}$ are state transition matrices, which are usually some state transformation operators;
$f(\cdot)$ is the objective function or fitness function, and
$y_{k+1}$ is the function value at $\bm x_{k+1}$.
\subsection{State transition operators}
Using state space transformation for reference, four special
state transition operators are designed to generate continuous solutions for an optimization problem.\\
(1) Rotation transformation
\begin{equation}
\bm x_{k+1}= \bm x_{k}+\alpha \frac{1}{n \|\bm x_{k}\|_{2}} R_{r} \bm x_{k},
\end{equation}
where $\alpha$ is a positive constant, called the rotation factor;
$R_{r}$ $\in$ $\mathbb{R}^{n\times n}$, is a random matrix with its entries being uniformly distributed random variables defined on the interval [-1, 1],
and $\|\cdot\|_{2}$ is the 2-norm of a vector. This rotation transformation
has the function of searching in a hypersphere with the maximal radius $\alpha$. \\
(2) Translation transformation\\
\begin{equation}
\bm x_{k+1} = \bm x_{k}+  \beta  R_{t}  \frac{\bm x_{k}- \bm x_{k-1}}{\|\bm x_{k}- \bm x_{k-1}\|_{2}},
\end{equation}
where $\beta$ is a positive constant, called the translation factor; $R_{t}$ $\in \mathbb{R}$ is a uniformly distributed random variable defined on the interval [0,1].
The translation transformation has the function of searching along a line from $x_{k-1}$ to $x_{k}$ at the starting point $x_{k}$ with the maximum length $\beta$.
\\
(3) Expansion transformation\\
\begin{equation}
\bm x_{k+1} = \bm x_{k}+  \gamma  R_{e} \bm x_{k},
\end{equation}
where $\gamma$ is a positive constant, called the expansion factor; $R_{e} \in \mathbb{R}^{n \times n}$ is a random diagonal
matrix with its entries obeying the Gaussian distribution. The expansion transformation
has the function of expanding the entries in $\bm x_{k}$ to the range of [-$\infty$, +$\infty$], searching in the whole space.\\
(4) Axesion transformation\\
\begin{equation}
\bm x_{k+1} = \bm x_{k}+  \delta  R_{a}  \bm x_{k}\\
\end{equation}
where $\delta$ is a positive constant, called the axesion factor; $R_{a}$ $\in \mathbb{R}^{n \times n}$ is a random diagonal matrix with its entries obeying the Gaussian distribution and only one random position having nonzero value. The axesion transformation is aiming to search along the axes, strengthening single dimensional search.
\subsection{Regular neighborhood and sampling}
For a given solution, a candidate solution is generated by using one of the aforementioned state transition operators.
Since the state transition matrix in each state transformation is random, the generated candidate solution is not unique.
Based on the same given point, it is not difficult to imagine that a regular neighborhood will be automatically formed
when using certain state transition operator.
In theory, the number of candidate solutions in the neighborhood is infinity;
as a result, it is impractical to enumerate all possible candidate solutions.

Since the entries in state transition matrix obey certain stochastic
distribution, for any given solution, the new candidate becomes a random vector and its corresponding solution
(the value of a random vector) can be regarded as a sample.
Considering that any two random state transition matrices in each state transformation
are independent, several times of state transformation (called the degree of search enforcement, \textit{SE} for short)
based on the same given solution are performed
for certain state transition operator, consisting of \textit{SE} samples.
It is not difficult to find that all of the \textit{SE} samples are independent, and they are
representatives of the neighborhood.
Taking the rotation transformation for example, a total number of
\textit{SE} samples are generated in pseudocode as follows
\begin{algorithmic}[1]
\For{$i\gets 1$, \textit{SE}}
\State $\mathrm{State}(:,i) \gets {\mathrm{Best}} +\alpha \frac{1}{n \| {\mathrm{Best}} \|_{2}} R_{r}  {\mathrm{Best}} $
\EndFor
\end{algorithmic}
where $Best$ is the incumbent best solution, and \textit{SE} samples are stored in
the matrix $\mathrm{State}$.

\subsection{An update strategy}
As mentioned above, based on the incumbent best solution, a total number of
\textit{SE} candidate solutions are generated.
A new best solution is selected from the candidate set by virtue of the fitness function, denoted as \emph{newBest}.
Then, an update strategy based on greedy criterion is used to update the incumbent best as shown below
\begin{equation}
Best =
\left\{ \begin{aligned}
&\mathrm{newBest},  \;\;\;\mathrm{if}\; f(\mathrm{newBest}) < f(\mathrm{Best}) \\
&\mathrm{Best},\;\;\;\;\;\;    \;\;\;\mathrm{otherwise}
\end{aligned} \right.
\end{equation}
\subsection{Algorithm procedure of the basic continuous STA}
With the state transformation operators, sampling technique and update strategy, the basic state transition algorithm can be described by the following pseudocode
\begin{algorithmic}[1]
\Repeat
    \If{$\alpha < \alpha_{\min}$}
    \State {$\alpha \gets \alpha_{\max}$}
    \EndIf
    \State {Best $\gets$ expansion(funfcn,Best,SE,$\beta$,$\gamma$)}
    \State {Best $\gets$ rotation(funfcn,Best,SE,$\alpha$,$\beta$)}
    \State {Best $\gets$ axesion(funfcn,Best,SE,$\beta$,$\delta$)}
    \State {$\alpha \gets \frac{\alpha}{\textit{fc}}$}
\Until{the specified termination criterion is met}
\end{algorithmic}

\indent As for detailed explanations, rotation$(\cdot)$ in above pseudocode is given for illustration purposes as follows
\begin{algorithmic}[1]
\State{oldBest $\gets$ Best}
\State{fBest $\gets$ feval(funfcn,oldBest)}
\State{State $\gets$ op\_rotate(Best,SE,$\alpha$)}
\State{[newBest,fnewBest] $\gets$ fitness(funfcn,State)}
\If{fnewBest $<$ fBest}
    \State{fBest $\gets$ fnewBest}
    \State{Best $\gets$ newBest}
    \State{State $\gets$ op\_translate(oldBest,newBest,SE,$\beta$)}
    \State{[newBest,fnewBest] $\gets$ fitness(funfcn,State)}
    \If{fnewBest $<$ fBest}
        \State{fBest $\gets$ fnewBest}
        \State{Best $\gets$ newBest}
    \EndIf
\EndIf
\end{algorithmic}

As shown in the above pseudocodes, the rotation factor $\alpha$ is decreasing periodically from a maximum
value $\alpha_{\max}$ to a minimum value $\alpha_{\min}$ in an
exponential way with base \textit{fc}, which is called lessening coefficient.
op\_rotate$(\cdot)$ and op\_translate$(\cdot)$ represent the implementations of proposed sampling technique for rotation and
translation operators, respectively, and fitness$(\cdot)$ represents the implementation of selecting the new best solution
from \textit{SE} samples. It should be noted that the translation operator is only executed
when a solution better than the incumbent best solution can be found in the \textit{SE} samples
from rotation, expansion or axesion transformation.
In the basic continuous STA, the parameter settings are given as follows:
$\alpha_{\max} = 1, \alpha_{\min} = 1e$-4, $\beta = 1, \gamma = 1, \delta = 1$, $\textit{SE} = 30, \textit{fc} = 2$.

When using the fitness$(\cdot)$ function, solutions in \emph{State} are projected into
$\Omega$ by using the following formula
\begin{equation}
x_i =
\left\{ \begin{aligned}
&u_i,  \;\;\; \mathrm{if}\; x_i > u_i \\
&l_i,  \;\;\; \mathrm{if}\; x_i < l_i\\
&x_i,  \;\;\; \mathrm{otherwise}
\end{aligned} \right.
\end{equation}
where $u_i$ and $l_i$ are the upper and lower bounds of $x_i$ respectively.

\section{A matlab implementation of the continuous STA}

\subsection{Installation}
The STA toolbox was developed under matlab 7.11.0 (R2010b),  which can be download via the following link
\url{http://www.mathworks.com/matlabcentral/fileexchange/52498-state-transition-algorithm}.

After unzipping the file ``basic\_STA.zip", you will get the following file list:
\begin{figure}[!htb]
  \centering
  \includegraphics[width=\hsize]{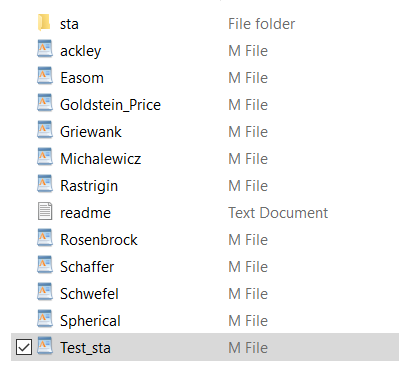}
  \caption{The unzipped basic\_STA file list}
  \label{fig1}
\end{figure}

There exist a file folder named ``sta" (it contains the core files of the STA toolbox),
and several .m files as well as a text document.

\subsection{The main file}
By running the main file ``Test\_sta.m", which is shown as below
\begin{verbatim}
clear all
clc
currentFolder = pwd;
addpath(genpath(currentFolder))
% parameter setting
warning('off')
SE =  30; % degree of search enforcement
Dim = 10;% dimension
Range = repmat([-5.12;5.12],1,Dim);
Iterations = 1e3;
tic
[Best,fBest,history] = STA(@Rastrigin,SE,
Dim,Range,Iterations);
toc
Best
fBest
semilogy(history)
xlabel('Iterations'),ylabel('Fitness(log)')
\end{verbatim}
you can
get the following similar results
\begin{figure}[!htb]
  \centering
  \includegraphics[width=\hsize]{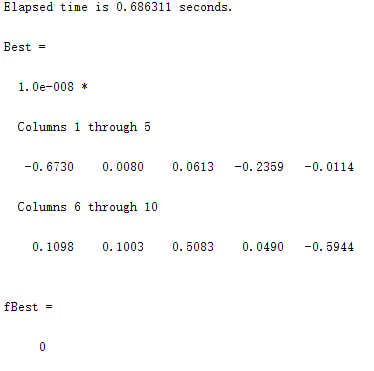}
  \includegraphics[width=\hsize]{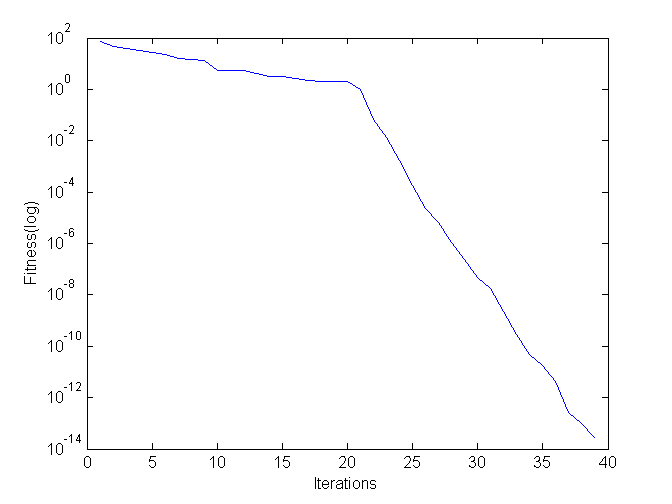}
  \caption{The results for the Rastrigin function (10D)}
  \label{fig2}
\end{figure}

If you want to optimize another benchmark function (the file list as shown in Fig. \ref{fig1}), for example, the
Griewank function with 15 dimension, you just need to do the following changes
\begin{verbatim}
Dim = 15;% dimension
Range = repmat([-600;600],1,Dim);%range
[Best,fBest,history] = STA(@Griewank,SE,
Dim,Range,Iterations);
\end{verbatim}

\begin{figure}[!htb]
  \centering
  \includegraphics[width=\hsize]{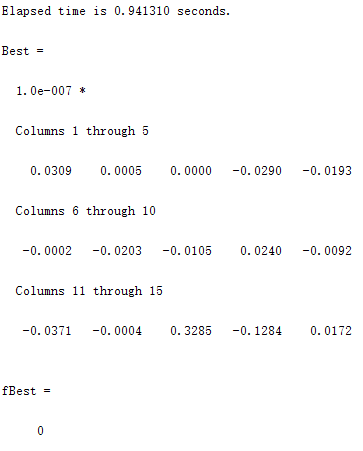}
  \includegraphics[width=\hsize]{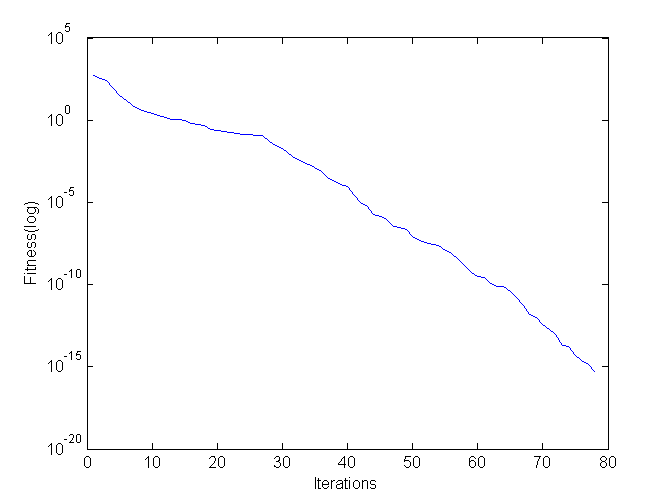}
  \caption{The results for the Griewank function (15D)}
  \label{fig3}
\end{figure}

After that, by running the main file, you can get the following similar results as shown in
Fig. \ref{fig3}.

If you want to optimize a user-defined function, for example,
the following minimization problem:
\begin{eqnarray*}
&&\min_{\bm x \in \mathbb{R}^3} \;\; (x_1 - 1)^2 + (x_2 - 2x_1)^2 + (x_3 - 3 x_2)^2 \nonumber \\
&&\mathrm{s.t.} \;\; -3 \leq x_1 \leq 3, -2 \leq x_2 \leq 2, -1 \leq x_3 \leq 1
\end{eqnarray*}

Firstly, you have to create a m-file (for example with a name `myfun') as follows
\begin{verbatim}
function y = myfun(x)
x1=x(:,1);
x2=x(:,2);
x3=x(:,3);
y=(x1-1).^2+(x2 - 2*x1).^2+(x3-3*x2).^2;
\end{verbatim}

Then, you need to do the following changes of the main file ``Test\_sta.m"
\begin{verbatim}
Dim = 3;% dimension
Range = [-3 -2 -1;3 2 1];%range
Iterations = 1e1;
tic
[Best,fBest,history] = STA(@myfun,SE,
Dim,Range,Iterations);
toc
Best
fBest
\end{verbatim}

\begin{figure}[!htb]
  \includegraphics[width=6cm]{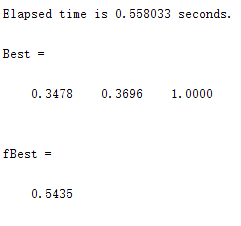}
  \includegraphics[width=\hsize]{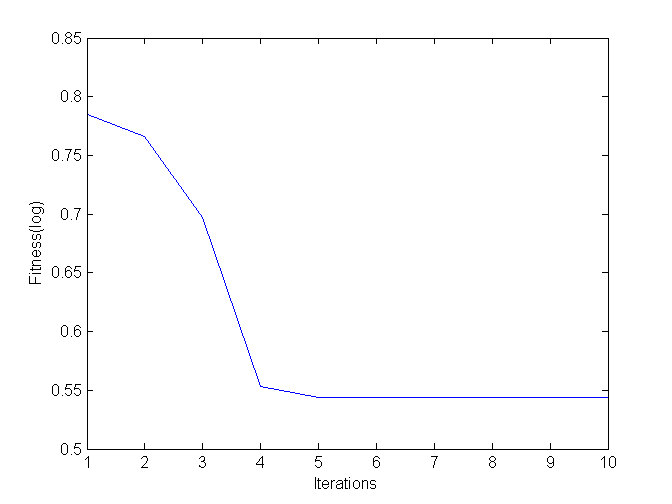}
  \caption{The results for the user-defined function (3D)}
  \label{fig4}
\end{figure}

After that, by running the main file, you can get the following similar results as shown in
Fig. \ref{fig4}.

\subsection{The core files of STA}
\begin{figure}[!htb]
  \includegraphics[width=\hsize]{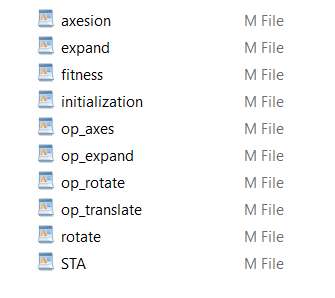}
  \caption{The file list in the file folder ``sta"}
  \label{fig5}
\end{figure}

The core files of STA are contained in the file folder ``sta", in which, the file list can be
seen in Fig. \ref{fig5}.

By opening the file STA.m, the main algorithm procedure of the basic continuous STA
is shown as below
\begin{verbatim}
function [Best,fBest,history]=STA(funfcn,
SE,Dim,Range,Iterations)
% parameter setting
alpha_max = 1;
alpha_min = 1e-4;
alpha = alpha_max;
beta = 1;
gamma = 1;
delta = 1;
fc = 2;
% initialization
State = initialization(SE,Dim,Range);
[Best,fBest] = fitness(funfcn,State);
% iterative process
for iter = 1:Iterations
    if alpha < alpha_min
        alpha = alpha_max;
    end
    [Best,fBest] = expand(funfcn,Best,
    SE,Range,beta,gamma);
    [Best,fBest] = rotate(funfcn,Best,
    SE,Range,alpha,beta);
    [Best,fBest] = axesion(funfcn,Best,
    SE,Range,beta,delta);
    history(iter) = fBest;
    alpha = alpha/fc;
end
\end{verbatim}

The function \textit{initialization} is to generate \textit{SE} uniformly random initial points within the range (lower and upper bounds),
and \textit{fitness} is to select the best candidate among them using greedy criterion. They are shown in matlab codes as below
\begin{verbatim}
function  State=initialization(SE,Dim,
Range)
Pop_Lb = Range(1,:);
Pop_Ub = Range(2,:);
State  = rand(SE,Dim).*repmat
(Pop_Ub-Pop_Lb,SE,1)+repmat(Pop_Lb,SE,1);

function [Best,fBest] = fitness(funfcn,
State) % calculate fitness
fState = feval(funfcn,State);
[fGBest, g] = min(fState);
fBest = fGBest;
Best = State(g,:);
\end{verbatim}

The functions \textit{op\_rotate}, \textit{op\_translate}, \textit{op\_expand} and \textit{op\_axes}
are the implementation of rotation, translation, expansion and axesion transformation respectively.
The corresponding matlab codes are listed as below
\begin{verbatim}
function y=op_rotate(Best,SE,alpha)
%rotation transformation
n = length(Best);
y = repmat(Best',1,SE) +
alpha*(1/n/(norm(Best)+eps))*
reshape(unifrnd(-1,1,SE*n,n)*Best',n,SE);
y = y';

function y=op_translate(oldBest,newBest,
SE,beta)
%translation transfomration
n = length(oldBest);
y = repmat(newBest',1,SE) +
beta/(norm(newBest-oldBest)+eps)
*reshape(kron(rand(SE,1),
(newBest-oldBest)'),n,SE);
y = y';

function y = op_expand(Best,SE,gamma)
%expansion transformation
n = length(Best);
y = repmat(Best',1,SE) + gamma*
(normrnd(0,1,n,SE).*repmat(Best',1,SE));
y = y';

function y = op_axes(Best,SE,delta)
% axesion transformation
n = length(Best);
A = zeros(n,SE);
index = randint(1,SE,[1,n]);
A(n*(0:SE-1)+index) = 1;
y = repmat(Best',1,SE) + delta*
normrnd(0,1,n,SE).*A.*repmat(Best',1,SE);
y = y';
\end{verbatim}

Each of these functions aims to generate \textit{SE}
random points with special property based on incumbent best solution.

The functions \textit{expansion}, \textit{rotation}, \textit{axesion} are to generate random points and update incumbent best
solution by using corresponding state transition operators. Taking the \textit{expansion} function
for example, it is shown as below
\begin{verbatim}
function [Best,fBest] = expand(funfcn,
oldBest,SE,Range,beta,gamma)
Pop_Lb=repmat(Range(1,:),SE,1);
Pop_Ub=repmat(Range(2,:),SE,1);
Best = oldBest;
fBest  = feval(funfcn,Best);
flag = 0;
State = op_expand(Best,SE,gamma);
changeRows = State > Pop_Ub;
State(find(changeRows))
= Pop_Ub(find(changeRows));
changeRows = State < Pop_Lb;
State(find(changeRows))
= Pop_Lb(find(changeRows));
[newBest,fGBest] = fitness(funfcn,State);
if fGBest < fBest
    fBest = fGBest;
    Best = newBest;
    flag = 1;
else
    flag = 0;
end

if flag ==1
    State =
    op_translate(oldBest,Best,SE,beta);
    changeRows = State > Pop_Ub;
    State(find(changeRows))
    = Pop_Ub(find(changeRows));
    changeRows = State < Pop_Lb;
    State(find(changeRows))
    = Pop_Lb(find(changeRows));
    [newBest,fGBest]
    = fitness(funfcn,State);
    if fGBest < fBest
        fBest = fGBest;
        Best = newBest;
    end
end
\end{verbatim}
and it can be found that the translation transformation is executed only when a better solution
is found by the expansion transformation (flag =1).

The matlab codes of \textit{rotation} function and \textit{axesion} function are very similar to that of
\textit{expansion} function. The only difference is that the subfunction \textit{op\_expand} is replaced by
\textit{op\_rotate} function or \textit{op\_axes} function.
For more details of the \textit{rotation} and \textit{axesion} functions, please refer to \cite{Zhou2015matlab}. 

\section{Conclusion}
In this paper, a matlab toolbox for the standard continuous STA is described, and 
several instances are given to show how to use the STA toolbox. 
The core functions in the STA toolbox are explained as well. 
An available link to download the STA toolbox is also provided for reference.

\end{document}